\documentclass[10pt,twosided]{article}
\usepackage[tbtags]{amsmath}
\usepackage{epsfig,amstext,amssymb,amsthm,latexsym}
\allowdisplaybreaks[4]
\usepackage{color}
\newcommand{\COM}[1]{}

\usepackage{amssymb,color}
\newcommand{\EE}[1]{\mathbb{E}\Bigl\{#1 \Bigr\}}

\definecolor{c20}{rgb}{0.,0.7,0.}
\definecolor{c30}{rgb}{0.,0.,1.}
\definecolor{c40}{rgb}{1,0.1,0.7}
\definecolor{c50}{rgb}{1,0,0}
\definecolor{c60}{rgb}{1,0.9,0.1}

\def\xx#1{\textcolor{c30}{#1}}
\def\xx#1{#1}

\def\aA#1{\textcolor{c30}{#1}}
\def\aA#1{#1}

\def\cE#1{\textcolor{c30}{#1}}
\def\cE#1{#1}
\def\ee#1{\textcolor{c30}{#1}}
\def\ee#1{#1}
\def\hh#1{\textcolor{c30}{#1}}
\def\hh#1{#1}
\def\cH#1{\textcolor{c30}{#1}}
\def\cH#1{#1}

\def\cW#1{{\textcolor{c40}{#1}}}
\def\cW#1{#1}
\def\cw#1{{\textcolor{c40}{#1}}}
\def\cw#1{#1}
\def\cz#1{{\textcolor{c40}{#1}}}
\def\cz#1{#1}
\def\aE#1{\textcolor{c30}{#1}}
\def\aH#1{\textcolor{c30}{#1}}
\def\aH#1{#1}
\def\aE#1{#1}
\def\cZ#1{{\textcolor{c40}{#1}}}
\def\cZ#1{#1}
\def\bE#1{\textcolor{c30}{#1}}
\def\bE#1{#1}
\def\CC#1{{\textcolor{c40}{#1}}}
\def\CC#1{#1}

\def\bC#1{{\textcolor{c40}{#1}}}
\def\bC#1{#1}
\def\CCC#1{{\textcolor{c40}{#1}}}
\def\CCC#1{#1}
\def\CCW#1{{\textcolor{c40}{#1}}}
\def\CCW#1{#1}

\usepackage{amssymb,color}
\def\lne{\ln n}

\newcommand{\ABs}[1]{ \biggl \lvert #1 \biggr \rvert}

\newcommand{\pk}[1]{\mathbb{P} \left( #1 \right) }

\newcommand{\EXP}[1]{\exp \left( #1 \right) }

\newcommand{\R}{\mathbb{R}}

\newcommand{\inr}{\in \R}

\newcommand{\limit}[1]{\lim_{#1 \to   \infty}}

\newcommand{\BQN}{\begin{eqnarray}}
\newcommand{\EQN}{\end{eqnarray}}
\newcommand{\BQNY}{\begin{eqnarray*}}
\newcommand{\EQNY}{\end{eqnarray*}}

\newcommand{\BS}{\begin{sat}}
\newcommand{\ES}{\end{sat}}
\newcommand{\BT}{\begin{theo}}
\newcommand{\ET}{\end{theo}}
\newcommand{\BK}{\begin{korr}}
\newcommand{\EK}{\end{korr}}

\newcommand{\BD}{\begin{de}}
\newcommand{\ED}{\end{de}}

\newcommand{\BIT}{\begin{itemize}}
\newcommand{\EIT}{\end{itemize}}
\newcommand{\BDI}{\begin{description}}
\newcommand{\EDI}{\end{description}}

\newcommand{\BRM}{\begin{remarks}}
\newcommand{\ERM}{\end{remarks}}

\newcommand{\BEL}{\begin{lem}}
\newcommand{\EEL}{\end{lem}}

\numberwithin{equation}{section}
\newtheorem{theo}{Theorem}[section]
\newtheorem{sat}[theo]{Proposition}
\newtheorem{de}[theo]{Definition}
\newtheorem{lem}{Lemma}[section]

\newtheorem{korr}[theo]{Corollary}

\newtheorem{remarks}[theo]{Remarks}

\newcommand{\QED}{\hfill $\Box$}

\newcommand{\netheo}[1]{{Theorem \ref{#1}}}

\def\P{\operatorname*{P}}
\def\Var{\operatorname*{Var}}
\def\Cov{\operatorname*{Cov}}

\def\I{\operatorname*{\mathbb{I}}}

\begin{document}

\centerline{\bf \large Limit Laws for Extremes of Dependent Stationary Gaussian Arrays}

  \vskip 1.5 cm
\centerline{Enkelejd Hashorva and Zhichao Weng }

  \vskip .1 cm
\centerline{Department of Actuarial Science, University of Lausanne, Switzerland}

  \vskip .5 cm

  \vskip .1 cm
{\bf Abstract:} In this paper we show that the componentwise maxima of
weakly dependent bivariate stationary Gaussian triangular arrays converge in distribution after normalisation to H\"usler-Reiss distribution. Under a strong dependence assumption, we prove that the limit distribution of the maxima is
a mixture of a bivariate Gaussian distribution and H\"usler-Reiss distribution.
Another finding of our paper is that the componentwise maxima and componentwise minima remain asymptotically independent even in the
settings of H\"usler and Reiss (1989) allowing further for weak dependence.
Further we derive an almost sure limit theorem under the Berman condition for the components of the triangular array.



\def\IF{\infty}

\section{Introduction and Main Result}
An important multivariate distribution of extreme value theory is the so-called H\"usler-Reiss distribution, which in a bivariate setting is given by \begin{eqnarray}\label{eq2.3}
H_{\lambda}(x,y)=\exp\left(-\Phi\left(\lambda+\frac{x-y}{2\lambda}\right)\exp(-y)-
\Phi\left(\lambda+\frac{y-x}{2\lambda}\right)\exp(-x)\right), \quad x,y\inr,
\end{eqnarray}
with $\Phi(\cdot)$ the univariate standard Gaussian distribution and $\lambda\in [0,\IF]$ a parameter. When $\lambda=0$ we have in fact
for any $x,y \inr$ that $H_0(x,y)= \min (\Lambda(x), \Lambda(y))$ and in the other extreme \aA{case} $\lambda=\IF$ the distribution function
$H_\IF$ is a product distribution with Gumbel marginals $\Lambda(x)=\exp(-\exp(-x)),x\inr$.
It is clear that for any $\lambda\in [0,\IF]$
the marginals of $H_\lambda$ are Gumbel distributions.\\
A striking property of $H_\lambda$ is that it is \cH{a} max-stable distribution function
(see Resnick (1987) for definition and main properties), and moreover it is a natural model for extremes of Gaussian triangular arrays,
as shown first in H\"usler and Reiss (1989), see also Falk et al. (2010). This fact is very important for statistical applications concerned with models for
extremes of dependent risks. The parameter $\lambda$ has a nice representation and comes naturally in the setup of Gaussian triangular arrays.
Roughly speaking, if $\rho(n)\in (-1,1)$ is the correlation coefficient of a bivariate triangular array, then under the H\"usler-Reiss condition
\begin{eqnarray}\label{eq1.1}
\cE{\lim_{n\to \IF} }(1-\rho(n)) \lne  = \lambda^2 \in [0,\infty]
\end{eqnarray}
the distribution function $H_\lambda$ appears as the limiting distribution of the \CCC{normalized} maxima.

In their seminal paper H\"usler and Reiss provided the full multivariate result where the condition \eqref{eq1.1} is assumed for any bivariate pair
of a multivariate Gaussian array.\\
 Recently, \aE{the research} interest on H\"usler-Reiss distribution has  grown significantly mainly due to the fact that
not only Gaussian, but chi-square, elliptical triangular arrays, and some more general models  have componentwise maxima attracted by that distribution (see Hashorva (2005), Frick and Reiss (2010), \aA{Hashorva (2008,2013)}, Hashorva et
al.\ (2012)). \aE{Additionally, new interesting applications have been proposed in
 Manjunath et al. (2012), and further crucial links with Brown-Resnick processes have been recently discovered, see the seminal contribution
Brown and Resncik (1977) for the first paper where the bivariate H\"usler-Reiss distribution appears, and see the recent deep contributions Kabluchko et al.\ (2009), Kabluchko (2011), Oesting et al.\ (2012), Engelke et al. (2014)
for the aforementioned links. We mention in passing that an important contribution somewhat related to the topic of our paper, but not to our techniques and proofs is O'Brien (1987).} \\
So far, \aE{primarily due to technical difficulties,} the available results in the literature \aE{are} concerned with the asymptotic behaviour of maxima of triangular arrays where the limit
distribution is H\"usler-Reiss \cz{distribution, under the assumption that} the components of triangular arrays are independent.

This contribution is the first attempt to allow for dependence among the components of triangular array, remaining in a Gaussian framework.
Specifically, we deal with
 $\bigl\{\mathbf{X}_{n,k}=\bigl(X_{n,k}^{(1)},X_{n,k}^{(2)}\bigr), 1\leq k\leq n, n\geq 1\bigr\}$
a triangular array of bivariate Gaussian random vectors with zero-mean, unit-variance
and correlation given by
\begin{eqnarray*}
corr\left(X_{n,k}^{(1)},X_{n,k}^{(2)}\right)=\rho_0(n),   \qquad
corr\left(X_{n,k}^{(i)},X_{n,l}^{(j)}\right)=\rho_{ij}(k,l,n)=\rho_{ij}(|k-l|,n),
\end{eqnarray*}
\aE{where $1\leq k\neq l\leq n$ and $i,j \in \{1,2\}.$}\\
For notational simplicity, in the sequel define the partial maxima
$$\left(M_{n,m}^{(1)},M_{n,m}^{(2)}\right):=\left(\max_{m+1\leq k\leq n}X_{n,k}^{(1)},\max_{m+1\leq k\leq n}X_{n,k}^{(2)}\right),
\qquad \left(M_{n}^{(1)},M_{n}^{(2)}\right):=\left(M_{n,0}^{(1)},M_{n,0}^{(2)}\right).$$
Motivated by Berman condition (see Berman (1964), or Berman (1992))
 we consider the following weak dependence condition adapted for the triangular array setup of our paper:

\COM{
{\bf Assumption A1}: Suppose that there exist positive integers $k_n<n, n\ge 1$ such that 
\begin{eqnarray}
\lim_{n\to \IF} k_n =\infty, \quad
\label{eq2.1}
\CC{\lim_{n\to \IF} \sup_{k_n\leq k < n \atop i,j \in \{1,2\}} \rho_{ij}(k,n) \ln  s =0.}
\end{eqnarray}
}

{\bf Assumption A1}: \bC{Suppose that $\sigma:=\max_{ 1\leq k < n,n\ge 2
\atop i,j \in \{1,2\} }|\rho_{ij}(k,n)|<1$. Let $  \alpha \in (0,
\frac{1-\sigma}{1+\sigma})$ and $I_n:=[n^{\alpha}]$, and further
assume that}
\begin{eqnarray}
\label{eq2.1}
\bC{\lim_{n\to \IF} \max_{I_n \leq k < n \atop i,j \in \{1,2\}} |\rho_{ij}(k,n)| \ln  n =0.}
\end{eqnarray}

Our first result below shows that the Berman condition does not change the limit distribution of the componentwise maxima, thus its asymptotic behaviour
is the same as in the iid setup of H\"usler and Reiss (1989). 

\BT\label{th2.1}
Let $\left\{\mathbf{X}_{n,k}, 1\leq k\leq n, n\geq 1\right\}$ be a bivariate Gaussian triangular array
satisfying Assumption A1. If \eqref{eq1.1} holds for $\rho_0(n)$ with $\lambda\in [0,\IF]$,
\COM{and further
\begin{eqnarray}\label{eq1.2}
\sup_{1\leq k\neq l\leq n, i,j\cE{\in} \{1,2\}}|\rho_{ij}(|k-l|,n)|< 1
\end{eqnarray}
holds,}
\cZ{then we have} 
\begin{eqnarray}\label{main:A}
\lim_{n \to \infty}\sup_{x,y\inr}
\ABs{\P\left(M_n^{(1)}\leq \ee{u_n(x)}, M_n^{(2)}\leq \ee{u_n(y)}\right) -H_{\lambda}(x,y)}=0, \quad
\end{eqnarray}
where $u_n(s)=a_ns+b_n,s\inr$ \bC{ with normalized constants $a_{n}$
and $b_{n}$ given by}
\begin{eqnarray}\label{eq2.2}
\cW{a_n=\frac{1}{\sqrt{2  \lne }}} \qquad \mbox{and} \qquad
b_n= \cE{\sqrt{2 \lne }  -
\frac{1}{2\sqrt{2  \lne } } \ln (4 \pi  \lne  )}.
\end{eqnarray}
\ET

A natural \ee{relaxation} of the weak dependence assumption is to allow the limit in \eqref{eq2.1} to be positive which we formulate \xx{below} as
\xx{our} second main
assumption, namely:
\COM{
{\bf Assumption A2}: Suppose that there exist positive integers $k_n<n, n\ge 1$ such that
\begin{eqnarray}
\lim_{n\to \IF} k_n =\infty, \quad
\label{eqST}
\CC{\lim_{n\to \IF} \sup_{k_n\leq s < n}\rho_{ij}(s,n) \ln  s =\tau_{ij}\in (0,\IF)}
\end{eqnarray}
for $i,j\in \{1,2\}$.}

{\bf Assumption A2}: Let $\tau_{ij}\in (0,\IF)$ be constants for $i,j\in \{1,2\}$. Suppose that $\delta_{ij}:=\max_{ 1\leq k < n, n\ge 2}|\rho_{ij}(k,n)|<1$, \bC{set} $\varpi_{ij} \in (0,
\frac{1-\delta_{ij}}{1+\delta_{ij}})$ and
$K_{n,ij}:=[n^{\varpi_{ij}}]$, \bC{and further
assume that}
\begin{eqnarray}
\label{eqST}
\lim_{n\to \IF} \max_{K_{n,ij}\leq k< n} |\rho_{ij}(k,n) \ln  k -\tau_{ij}|=0.
\end{eqnarray}

\ee{As shown in our second result below, in} the case of strong dependence, the limiting distribution of the joint maxima \cH{is given by a Gaussian \cZ{distribution} mixture of the H\"usler-Reiss} \cZ{distribution}.
We set next
$ \widetilde \tau:= \tau_{12} -\frac{1}{2}(\tau_{11}+\tau_{22})$ and $\widetilde{\lambda}:=\sqrt{\lambda^2+\widetilde \tau}$.

\BT\label{th2.3} Under the assumptions of \netheo{th2.1}, if further Assumption A2 holds,
\CCC{and assume that} $\tau_{12}\le \sqrt{\tau_{11}\tau_{22}}$ and $\lambda^2 \cz{\geq}- \widetilde \tau$,
then we have
\begin{eqnarray}
\lim_{n \to \infty}\sup_{x,y\inr} \ABs{\P\left(M_n^{(1)}\leq \ee{u_n(x)} , M_n^{(2)}\leq \ee{u_n(y)}\right)
-\EE{H_{\widetilde{\lambda}}\left(x+\tau_{11}-\sqrt{2\tau_{11}}Z,y+\tau_{22}-\sqrt{2\tau_{22}}W\right)}}=0, \quad
\end{eqnarray}
with 
$(Z,W)$ a standard bivariate Gaussian random vector with correlation $\tau_{12}/\sqrt{\tau_{11}\tau_{22}}$.
\ET

\ee{\CCW{The rest of the paper is organized as follows}: In Section 2 we discuss the novelty and \cH{the} importance of our results as well as connections with available contributions in the literature. \bE{Further we provide two important extensions:  first we show that the componentwise minima is asymptotically independent of the componentwise maxima. This result is well-known for the case of iid bivariate Gaussian sequences. This article is the first to \CCW{consider the asymptotic behavior of maxima and minima of stationary bivariate Gaussian arrays under the H\"{u}sler and Reiss condition and weak dependent condition.} 
Our second result in Section 2}
derives \cH{the} almost sure \bE{limit theorem} for the case of weak dependence. Proofs and further results are relegated to Section 3.}

\section{Discussion and Extensions}
As mentioned in the Introduction, all the contributions so far have considered only independent triangular arrays. Our motivation to allow
dependence comes from practical situations, where due to the presence of some random inflation/deflation or measurement errors (which are always present) the independence of the components of triangular arrays is not an adequate assumption. This is the first contribution in this direction;
we have treated the classical Gaussian setup since the dependence in more general models is very difficult to deal with.
Both findings displayed above are of interest: in case of weakly dependent bivariate stationary Gaussian triangular arrays the limiting distribution is H\"usler-Reiss distribution, which is identical to the iid case. Our second result shows that this is no longer the case for strongly dependent stationary Gaussian arrays. \\
An immediate consequence of \netheo{th2.3} is that the univariate maxima also converges after appropriate \CCC{normalization}, i.e.,
\begin{eqnarray}
\lim_{n\to \IF} \P\left(M_n^{(1)}\leq \ee{u_n(x)}\right) = \EE{\Lambda(x+\tau_{11}-\sqrt{2\tau_{11}}Z)}, \quad x\inr,
\end{eqnarray}
which is already derived in   Corollary 6.5.2 in Leadbetter  et al.\ (1983).\\
It can be easily seen that our results \aA{hold} for multivariate setup and not only for the bivariate setup; we refrain ourself to the bivariate setup for ease of presentation.\\
 There are different (interesting) possibilities to continue the investigation under dependence.
For instance, as in Hashorva (2011), one direction is to investigate if the convergence in \eqref{main:A} can be stated in a stronger form as convergence of corresponding density functions.\\
\bE{Our first result below concerns the joint asymptotic behaviour of the sample maxima and sample minima. For the iid setup, it is well-known that for multivariate Gaussian random sequences, \CCC{the} sample maxima and sample minima are asymptotically independent \CCW{(see Davis (1979))}. In the framework of triangular arrays suggested by H\"usler and Reiss (1989) no investigation in this direction has been done. Our first result below shows that the asymptotic independence is preserved even in the case of weakly dependent Gaussian arrays.
\BT\label{th1.2} Under the assumptions of \netheo{th2.1} we have
\BQN
\lim_{n \to \infty}\sup_{x_1,x_2,y_1,y_2\inr}
\Biggl \lvert \pk{-u_n(y_1)< m_n^{(1)}\leq M_n^{(1)}\leq u_n(x_1),-u_n(y_2)< m_n^{(2)}\leq M_n^{(2)}\leq u_n(x_2)}
\notag\\
&&\hspace{-11 cm} - H_{\lambda}(x_1,x_2)H_{\lambda}(y_1,y_2)\Biggr \lvert =0,
\EQN
with $\left(m_{n}^{(1)},m_{n}^{(2)}\right):=\left(\min_{1\leq k\leq n}X_{n,k}^{(1)},\min_{1\leq k\leq n}X_{n,k}^{(2)}\right)$ the componentwise sample minima.
\ET}
\COM{
\CCC{\BK
Under the assumptions of Theorem \ref{th1.2}, we have
\BQN
\lim_{n \to \infty}\sup_{x_1,x_2\inr}
\Biggl \lvert \pk{a_n^{-1}\left(\max_{1\leq s \leq n}|X_{n,s}^{(i)}|-b_n\right)\leq x_i+\ln 2, 1\leq i \leq 2}
 - H_{\lambda}(x_1,x_2)\Biggr \lvert =0.
\EQN
\EK
}
}
A different direction which we pursue below  is to analyse whether the convergence in \eqref{main:A}
can be strengthen to almost sure limit convergence; \cH{some related results in this direction are obtained in \CCC{Cheng et al. (1998),
 Fahrner and Stadtm\"{u}ller (1998), Cs\'{a}ki and Gonchigdanzan (2002),} Tan et al. (2007), Tan and Peng (2009),
Peng et al. (2010,2011),} Tan and Wang (2011,2012),  Weng et al. (2012).

\BT\label{th2.2} Under the assumptions and notation of \netheo{th2.1}, suppose further that the elements from different rows
satisfy
\BQN\label{add eq1}
corr\left(X^{(i)}_{mk}, X^{(j)}_{nl}\right)=\gamma_{ij}(k,l,m,n)=\gamma_{ij}(|k-l|+1,m,n) 
\EQN
with $1\leq k\leq m, 1\leq l \leq n$, and all $m < n, i,j\in\{1,2\}$.
Let $\theta:=\max_{1\leq k\leq n, n\ge 2, i,j\in \{1,2\}}|\gamma_{ij}(k,m,n)|<1$, $\beta \in (0, \frac{1-\theta}{1+\theta})$ and
$J_n:=[n^{\beta}]$. If \aA{for some $\epsilon>0$}
\BQN \label{eq:k_n}
\bC{\max_{I_n\leq s <n }\rho_{ij}(s,n)\ln s(\ln\ln s)^{1+\epsilon}=O(1)}
\EQN
and
\BQN\label{eq:l_n}
\bC{\max_{J_n\leq t \leq n}\gamma_{ij}(t,m,n)\ln t(\ln\ln t)^{1+\epsilon}=O(1)}
\EQN
hold as $n \to \infty$, then for any $(x,y) \in \ee{\R}^2$
\BQN\label{eq2.4}
\lim_{n \to \infty}\frac{1}{ \lne }\sum^n_{k=1}\frac{1}{k}\I\left(M^{(1)}_k\leq u_k(x), M^{(2)}_k\leq u_k(y)\right)=H_{\lambda}(x,y)
\EQN
holds almost surely. Furthermore, for any $x_1,x_2,y_1,y_2 \in \R$, almost surely
\BQN\label{eq2.7}
&&\lim_{n \to \infty}\frac{1}{ \lne }\sum^n_{k=1}\frac{1}{k}
\I\left(-u_k(y_1)< m_k^{(1)}\leq M_k^{(1)}\leq u_k(x_1),-u_k(y_2)< m_k^{(2)}\leq M_k^{(2)}\leq u_k(x_2)\right)\nonumber\\
&=&H_{\lambda}(x_1,x_2)H_{\lambda}(y_1,y_2).
\EQN
\ET
\COM{
\BK
Under the assumptions of Theorem \ref{th2.2}, for any $x_1,x_2 \in \R$ we have
\BQNY
\lim_{n \to \infty}\frac{1}{ \lne }\sum^n_{k=1}\frac{1}{k}
\I\left(a_k^{-1}\left(\max_{1\leq s \leq k}|X_{k,s}^{(i)}|-b_k\right)\leq x_i +\ln 2, 1\leq i \leq 2 \right)
=H_{\lambda}(x_1,x_2).
\EQNY
\EK
}

Strengthening of the result of \netheo{th2.3} to almost sure limit theorem requires significantly more efforts and additional technical conditions,
therefore we shall not address that point here.

\section{Further Results and Proofs}
In this section we present the \bC{proofs} of the main results.
Since those proofs depend on some results which are of some
independent interest, we formulate several lemmas.

\BEL\label{le3.1}
Under Assumption A1 \ee{for $i,j \in \{1,2\}$ and any $x,y \inr$}
we have
\BQN\label{L1}
\lim_{n \to \infty}n \sum_{k=1}^{n-1} |\rho_{ij}(k,n)| \EXP{-\frac{\omega^2_n}{1+|\rho_{ij}(k,n)|}}=0, \quad
\EQN
where   $\omega_n:=\min(|u_n(x)|, |u_n(y)|)$.
If Assumption A2 holds, then for $i,j \in \{1,2\}$
\BQN\label{L2}
\lim_{n \to \infty}n\sum^{n-1}_{k=1}|\rho_{ij}(k,n)-\tau_{ij}(n)|\exp\left(-\frac{\omega_n^2}{1+\varrho_{ij}(k,n)}\right)=0,
\EQN
where
$\tau_{ij}(n):=\tau_{ij}/\ln n$ and $\varrho_{ij}(k,n)=\max\{|\rho_{ij}(k,n)|,\tau_{ij}(n)\}$.
\EEL

\textbf{Proof.}
\ee{Our proof is similar to that of Lemma 4.3.2 and Lemma 6.4.1 in Leadbetter  et al.\ (1983)}. For notational simplicity, we omit the index and write below simply $\rho(k,n)$ instead of $\rho_{ij}(k,n)$,
and similarly we write $\tau, \tau(n), \varrho(k,n), \delta, \varpi, K_n$ for $\tau_{ij}, \tau_{ij}(n),
\varrho_{ij}(k,n), \delta_{ij}, \varpi_{ij}, K_{n,ij}$ respectively.
First note that \cE{in view of}  \eqref{eq2.2}
\BQNY
\lim_{n\to\infty}\frac{n}{\sqrt{2\pi} \ee{b_n}}\EXP{-\frac{u^2_n(s)}{2}}=\exp(-s), \, s\inr, \text{ and }
\lim_{n\to\infty}\frac{\omega_{n}}{\sqrt{2 \lne }}=\lim_{n\to\infty}\frac{b_{n}}{\sqrt{2 \lne }}=1.
\EQNY
\aH{Next, we write}
\BQNY
n\sum_{k=1}^{n-1}|\rho(k,n)|\exp\left(-\frac{\omega_{n}^2}{1+|\rho(k,n)|}\right)
=n\left(\sum_{k=1}^{I_n}+\sum_{k=I_n+1}^{n-1}\right)
|\rho(k,n)|\exp\left(-\frac{\omega_{n}^2}{1+|\rho(k,n)|}\right)=:S_{n1}+ S_{n2}.
\EQNY
\aH{with $I_n$ given by condition formulated in Assumption A1.}
By the choice of $\alpha$ and $\sigma$ in Assumption A1, for some positive  constants $c_1,c_2$
\BQNY
S_{n1} \le nn^{\alpha}\EXP{-\frac{\omega_{n}^2}{1+\sigma}} &
=&n^{1+\alpha}\left(\EXP{-\frac{\omega_{n}^2}{2}}\right)^{\frac{2}{1+\sigma}}\nonumber\\
&\leq&c_1 n^{1+\alpha}\left(\frac{\omega_{n}}{n}\right)^{\frac{2}{1+\sigma}}
\leq c_2 n^{1+\alpha-\frac{2}{1+\sigma}}( \lne )^{\frac{1}{1+\sigma}} \nonumber\\
&\to & 0,\quad n\to \IF.
\EQNY
\aE{Next,} define $\sigma(l,n):=\max_{l\leq k <n }|\rho(k,n)|<1$, \aA{then}  \CC{by \eqref{eq2.1}}
\BQNY
\lim_{n \to \infty}\sigma(I_n,n)\omega_{n}^2=\lim_{n \to \infty}2\sigma(I_n,n) \lne=0, 
\EQNY
\aA{and hence} we have 
\BQNY
S_{n2}
&\leq&n\sigma(I_n,n)\EXP{-\omega_{n}^2}\sum^{n-1}_{k=I_n+1}
\EXP{\frac{\omega_{n}^2|\rho(k,n)|}{1+|\rho(k,n)|}}\\
&\leq&n^2\sigma(I_n,n)\EXP{-\omega_{n}^2}\EXP{\sigma(I_n,n)\omega_{n}^2}\\
&\leq&O \Bigl(\sigma(I_n,n)\omega_{n}^2\EXP{\sigma(I_n,n)\omega_{n}^2}\Bigr)\\
&\to& 0, \quad n\to \IF,
\EQNY
\CCC{thus} \eqref{L1} follows. The following constant
$$
\delta(l,n):=\max_{l\leq k <n}\varrho(k,n)<1
$$
\aH{plays an important role for proving \eqref{L2}. We split the sum in \eqref{L2} into two terms, the first consist of summation over  $1\leq k\leq K_n$ and the second term is the sum over $K_n< k \CCC{<} n$. As in the proof of $S_{n1}$ above, for some positive constants $c,c_1,c_2$ we have
}
\BQNY
&&n\sum_{k=1}^{K_n}
|\rho(k,n)-\tau(n)|\exp\left(-\frac{\omega_{n}^2}{1+\varrho(k,n)}\right)\\
&\le& cnn^{\varpi}\EXP{-\frac{\omega_{n}^2}{1+\max(\delta,\tau(n))}} \\&
=&cn^{1+\varpi}\left(\EXP{-\frac{\omega_{n}^2}{2}}\right)^{\frac{2}{1+\max(\delta,\tau(n))}}\nonumber\\
&\leq&c_1 n^{1+\varpi}\left(\frac{\omega_{n}}{n}\right)^{\frac{2}{1+\max(\delta,\tau(n))}}
\leq c_2 n^{\varpi-\frac{1-\max(\delta,\tau(n))}{1+\max(\delta,\tau(n))}}( \lne )^{\frac{1}{1+\max(\delta,\tau(n))}} \nonumber\\
&\to & 0,\quad n\to \IF,
\EQNY
since $\tau(n)\to 0$ and $0< \varpi<\frac{1-\delta}{1+\delta}$.
For the second term, note that
\BQNY
&&n\sum^{n-1}_{k=K_n+1}|\rho(k,n)-\tau(n)|\exp\left(-\frac{\omega_n^2}{1+\varrho(k,n)}\right)\\
&\leq&n\exp\left(-\frac{\omega_n^2}{1+\delta(K_n,n)}\right)\sum^{n-1}_{k=K_n+1}|\rho(k,n)-\tau(n)|\\
&=&\frac{n^2}{ \lne }\exp\left(-\frac{\omega_n^2}{1+\delta(K_n,n)}\right)
\frac{ \lne }{ n}\sum^{n-1}_{k=K_n+1}|\rho(k,n)-\tau(n)|.
\EQNY
Since $\lim_{n \to \infty}\rho(K_n,n) \ln  K_n=\tau$,
there exists a constant $L>0$ such that $\rho(K_n,n) \ln  K_n\leq L$,
also $\delta(K_n,n) \ln  K_n \leq L$. \ee{Consequently, for some positive constants $c_1,c_2$} 
\BQNY
\frac{n^2}{ \lne }\exp\left(-\frac{\omega_n^2}{1+\delta(K_n,n)}\right)
&\leq&\frac{n^2}{ \lne }\exp\left(-\frac{\omega_n^2}{1+L/ \lne ^{\varpi}}\right)\\
&\leq&c_1 \frac{n^2}{ \lne }\left(\frac{\omega_n}{n}\right)^{2/(1+L/ \lne ^{\varpi})}\\
&\leq&c_2n^{2L/(L+ \lne ^{\varpi})}( \lne )^{-L/(L+ \lne ^{\varpi})}\\
&=&O(1), \quad n \to \infty.
\EQNY
Further
\BQNY
\frac{ \lne }{ n}\sum^{n-1}_{k=K_n+1}|\rho(k,n)-\tau(n)| &\leq&\frac{ \lne }{ n}\sum^{n-1}_{k=K_n+1}\left|\rho(k,n)-\frac{1}{ \ln k }\tau\right|
+\tau\frac{ \lne }{ n}\sum^{n-1}_{k=K_n+1}\left|\frac{1}{ \ln k }-\frac{1}{ \lne }\right|\\
&\leq&\frac{1}{ n\varpi}\sum^{n-1}_{k=K_n+1}\left|\rho(k,n) \ln k -\tau\right|
+\tau\frac{1}{ n}\sum^{n-1}_{k=K_n+1}\left|1-\frac{ \lne }{ \ln k }\right|\\
&=:&T_{n1}+T_{n2}.
\EQNY
By Assumption A2  $\lim_{n \to \IF} T_{n1}=0$, and further 
\BQNY
T_{n2} 
\leq\frac{\tau}{\varpi \lne }\sum^{n-1}_{k=K_n+1}\left| \ln  \frac{k}{n}\right|\frac{1}{n}
=O\Bigl(\frac{\tau}{\varpi \lne }\int^1_0| \ln  x| dx \Bigr),
\EQNY
thus the proof is complete.
\QED

\BEL\label{le3.2}
Under the conditions of Theorem \ref{th2.2} and the notation of Lemma \ref{le3.1},
for $i,j \in \{1,2\}$ and positive constants $c_1$, $c_2$
\BQN\label{eq3.1}
n\sum_{k=1}^{n-1}|\rho_{ij}(k,n)|\exp\left(-\frac{\omega_{n}^2}{1+|\rho_{ij}(k,n)|}\right)\leq c_1(\ln\ln n)^{-(1+\epsilon)},
\EQN
and
\BQN\label{eq3.2}
\max_{1\leq m < n}m\sum_{k=1}^{n}|\gamma_{ij}(k,m,n)|
\exp\left(-\frac{\omega_{m}^2+\omega_{n}^2}{2(1+|\gamma_{ij}(k,m,n)|)}\right)\leq c_2(\ln\ln n)^{-(1+\epsilon)}
\EQN
\cH{hold for some  $\epsilon>0$.}
\EEL

\textbf{Proof.}
For $m<n$, \cZ{using the notation of the proof of Lemma 3.1, and write simply $\gamma(k,m,n)$ instead of $\gamma_{ij}(k,m,n)$,}
further define the following constant
\BQNY
\theta(l,m,n):=\max_{l\leq k\leq n}|\gamma(k,m,n)|<1.
\EQNY
\COM{
and \hh{set} 
$J_n:=[n^{\beta}]$.
\hh{Clearly, we have}
$$\lim_{n \to \infty}I_n=\lim_{n \to\infty}J_n=\infty. $$}
By \eqref{eq:k_n} and $\eqref{eq:l_n}$, for some small $\varepsilon>0$ and some positive constants $c_1$, $c_2$ 
\cZ{and for all large $n$}
\BQNY
\sigma(I_n,n)\omega_{n}^2 &\leq& (1+\varepsilon)\frac{2}{\alpha}\sigma(I_n,n)\ln n^{\alpha}\\
&\leq& \bC{(1+\varepsilon)\frac{2}{\alpha}\max_{I_n\leq k< n}|\rho(k,n)|\ln k}\\
&\leq& \bC{c_{1}( \ln\ln n^{\alpha})^{-(1+\epsilon)} \sim c_1(
\ln\ln n)^{-(1+\epsilon)}}, \EQNY and for $1\leq m < n$ \BQNY
\theta(J_n,m,n)\omega_{m}\omega_{n}&\leq&
\bC{2(1+\varepsilon)\theta(J_n,m,n)(\ln m\ln n)^{1/2}}\\
&\leq& \bC{(1+\varepsilon)\frac{2}{\beta}\max_{ J_n\leq k \leq
n}|\gamma(k,m,n)|\ln k}\\
&\leq& \bC{c_{2}( \ln\ln n^{\beta})^{-(1+\epsilon)} \sim c_2( \ln\ln
n)^{-(1+\epsilon)}}.
\EQNY
\hh{Similarly}, \cH{for all $n$ large}
\BQNY \theta(J_n,m,n)\omega_{m}^2\leq
2(1+\varepsilon)\theta(J_n,m,n)\ln m \leq
(1+\varepsilon)\frac{2}{\beta}\theta(J_n,m,n)\ln n^{\beta} \leq c_2(
\ln\ln n)^{-(1+\epsilon)} \EQNY and \BQNY
\theta(J_n,m,n)\omega_{n}^2\leq 2(1+\varepsilon)\theta(J_n,m,n)\ln n
= (1+\varepsilon)\frac{2}{\beta}\theta(J_n,m,n)\ln n^{\beta} \leq
c_2( \ln\ln n)^{-(1+\epsilon)}. \EQNY
 \cZ{Combining the above inequalities and along the same lines of} the proof of Lemma 2.1 in Cs\'{a}ki and Gonchigdanzan (2002), the claim follows. \QED

\COM{
For the inequality \eqref{eq3.1}, split it into two parts. The first for $1\leq k\leq r_n$, it tends to zero as in Lemma \ref{le3.1}.
The second for $r_n< k <n$, note that
$$\lim_{n \to \infty}r_n=\infty$$
and
\BQNY
\sigma(r_n,n)\omega_{n}^2\leq \mathbb{C}\sigma(r_n,n)\ln n^{\tau}
\leq \mathbb{C}\sup_{k\geq r_n}\rho(k,n)\ln k
\leq \mathbb{C}( \ln\ln k)^{-(1+\epsilon)}
\leq \mathbb{C}( \ln\ln n^{\tau})^{-(1+\epsilon)}
\leq \mathbb{C}( \ln\ln n)^{-(1+\epsilon)}.
\EQNY
Hence we have
\BQNY
n\sum^n_{k=r_n+1}\rho(k,n)\EXP{-\frac{\omega_{n}^2}{1+\rho(k,n)}}
\leq\mathbb{C}\sigma(r_n,n)\omega_{n}^2\EXP{\sigma(r_n,n)\omega_{n}^2}
\leq \mathbb{C}( \ln\ln n)^{-(1+\epsilon)}.
\EQNY
For the inequality \eqref{eq3.2}, using similar arguments in Lemma 2.1 in Cs\'{a}ki and Gonchigdanzan (2002).
For $m<n$, $l< n,$ define next the following constants
$$\theta:=\sup_{1\leq n,1\leq k<n}\gamma(k,m,n)<1, \qquad \theta(l,m,n)=\sup_{k\geq l}\gamma(k,m,n)<1,\qquad
\xi =\frac{1-\theta}{2(1+\theta)}$$
and set $I_n=[n^{\xi}]$
\BQNY
m\sum_{k=1}^n\gamma(k,m,n)\exp\left(-\frac{\omega_{m}^2+\omega_{n}^2}{2(1+\gamma(k,m,n))}\right)
=m\left(\sum_{k=1}^{I_n}+\sum_{k=I_n+1}^n\right)
\gamma(k,m,n)\exp\left(-\frac{\omega_{m}^2+\omega_{n}^2}{2(1+\gamma(k,m,n))}\right)=:S_{n1}+S_{n2}.
\EQNY
By the choice of $\theta$ and $\xi$ we have
\BQNY
S_{n1}&\leq&mn^{\xi}\EXP{-\frac{\omega_{m}^2+\omega_{n}^2}{2(1+\theta)}}\\
&=&mn^{\xi}\left(\EXP{-\frac{\omega_{m}^2+\omega_{n}^2}{2}}\right)^{\frac{1}{1+\theta}}\\
&\leq&\mathbb{C}mn^{\xi}\left(\frac{\omega_{m}\omega_{n}}{mn}\right)^{\frac{1}{1+\theta}}\\
&\leq&\mathbb{C}m^{1-\frac{1}{1+\theta}}n^{\xi-\frac{1}{1+\theta}}(\ln m\ln n)^{\frac{1}{2(1+\theta)}}\\
&\leq&\mathbb{C}n^{1+\xi-\frac{2}{1+\theta}}(\ln n)^{\frac{1}{1+\theta}}\\
&\to & 0,  \qquad n \to \infty.
\EQNY
Further,
\BQNY
S_{n2}&\leq&m\theta(I_n,m,n)\EXP{-\omega_{m}^2-\omega_{n}^2}\sum^n_{k=I_n+1}
\EXP{\frac{(\omega_{m}^2+\omega_{n}^2)|\gamma(k,m,n)|}{1+|\gamma(k,m,n)|}}\\
&\leq&mn\theta(I_n,m,n)\EXP{-\omega_{m}^2-\omega_{n}^2}\EXP{\theta(I_n,m,n)(\omega_{m}^2+\omega_{n}^2)}\\
&\leq&\mathbb{C}\theta(I_n,m,n)\omega_{m}\omega_{n}\EXP{\theta(I_n,m,n)(\omega_{m}^2+\omega_{n}^2)}.
\EQNY
Notice that $\lim_{n \to\infty}I_n=\infty$,
\BQNY
\theta(I_n,m,n)\omega_{m}\omega_{n}&\leq& \mathbb{C}\theta(I_n,m,n)(\ln m\ln n)^{1/2}\\
&\leq& \mathbb{C}\theta(I_n,m,n)\ln n^{\xi}
\leq \mathbb{C}\sup_{k\geq I_n}\gamma(k,m,n)\ln k
\leq \mathbb{C}( \ln\ln k)^{-(1+\epsilon)}
\leq \mathbb{C}( \ln\ln n^{\xi})^{-(1+\epsilon)}\\
&\leq& \mathbb{C}( \ln\ln n)^{-(1+\epsilon)},
\EQNY
and for $1\leq m \leq n$,
\BQNY
\theta(I_n,m,n)\omega_{m}^2\leq \mathbb{C}\theta(I_n,m,n)\ln m \leq \mathbb{C}\theta(I_n,m,n)\ln n
\leq \mathbb{C}\theta(I_n,m,n)\ln n^{\xi}
\leq \mathbb{C}( \ln\ln n)^{-(1+\epsilon)}.
\EQNY
as $n \to \infty$. Thus the inequality \eqref{eq3.2} holds.
\QED
}

\textbf{Proof of Theorem \ref{th2.1}.}
Let $\bigl\{\bigl(\hat{X}_{n,k}^{(1)},\hat{X}_{n,k}^{(2)}\bigr), 1\leq k\leq n, n\geq 1 \bigr\}$
denote the associated iid triangular array of $\{\mathbf{X}_{n,k}\}$, i.e. the correlation satisfy
$corr\left(\hat{X}_{n,k}^{(1)},\hat{X}_{n,k}^{(2)}\right)=\rho_0(n)$ and
$corr\left(\hat{X}_{n,k}^{(i)},\hat{X}_{n,l}^{(j)}\right)=0$ for $1\leq k\neq l\leq n, i,j\in \{1,2\}$.
Since the condition \eqref{eq1.1} holds, using Theorem 1 in H\"{u}sler and Reiss (1989), we have
$$\lim_{n \to \infty}\sup_{x,y \in \R}\left|\P\left(\mathfrak{M}_n^{(1)}\leq u_n(x),\mathfrak{M}_n^{(2)}\leq u_n(y)\right)-H_{\lambda}(x,y)\right|=0,$$
where $\left(\mathfrak{M}_n^{(1)},\mathfrak{M}_n^{(2)}\right):=\left(\max_{1\leq k \leq n}\hat{X}_{n,k}^{(1)}, \max_{1\leq k \leq n}\hat{X}_{n,k}^{(2)}\right)$.
Hence, we only need to prove that
$$\lim_{n \to \infty}\sup_{x,y \in \R}\left|\P\left(M_n^{(1)}\leq u_n(x),M_n^{(2)}\leq u_n(y)\right)-
\P\left(\mathfrak{M}_n^{(1)}\leq u_n(x),\mathfrak{M}_n^{(2)}\leq u_n(y)\right)\right|=0$$
holds. By \aA{Berman's} Normal Comparison Lemma (\aA{see Piterbarg (1996) for generalised Berman inequality} and Corollary 2.1 in Li and Shao (2002)), for all $x,y \in \R$
we have
\begin{eqnarray*}
&&\left|\P\left(M_n^{(1)}\leq u_n(x), M_n^{(2)}\leq u_n(y)\right)-
\P\left(\mathfrak{M}_n^{(1)}\leq u_n(x),\mathfrak{M}_n^{(2)}\leq u_n(y)\right)\right|\\
\COM{
&\leq&\CC{\frac{1}{4}\sum_{1\leq k< l\leq n}|\Cov(X_{n,k}^{(1)},X_{n,l}^{(1)})|\exp\left(-\frac{u_n^2(x)}{1+|\Cov(X_{n,k}^{(1)},X_{n,l}^{(1)})|}\right)}\\
&&+\CC{
\frac{1}{4}\sum_{1\leq k\neq l\leq n}|\Cov(X_{n,k}^{(1)},X_{n,l}^{(2)})|\exp\left(-\frac{u_n^2(x)+u_n^2(y)}{2(1+|\Cov(X_{n,k}^{(1)},X_{n,l}^{(2)})|)}\right)}\\
&&\CC{+\frac{1}{4}\sum_{1\leq k< l\leq n}|\Cov(X_{n,k}^{(2)},X_{n,l}^{(2)})|\exp\left(-\frac{u_n^2(y)}{1+|\Cov(X_{n,k}^{(2)},X_{n,l}^{(2)})|}\right)}\\
}
&\leq&\frac{1}{4}n\sum_{k=1}^{n-1}|\rho_{11}(k,n)|\exp\left(-\frac{u_n^2(x)}{1+|\rho_{11}(k,n)|}\right)+
\frac{1}{2}n\sum_{k=1}^{n-1}|\rho_{12}(k,n)|\exp\left(-\frac{u_n^2(x)+u_n^2(y)}{2(1+|\rho_{12}(k,n)|)}\right)\\
&&+\frac{1}{4}n\sum_{k=1}^{n-1}|\rho_{22}(k,n)|\exp\left(-\frac{u_n^2(y)}{1+|\rho_{22}(k,n)|}\right).
\end{eqnarray*}
\aA{Consequently,} in view of \eqref{L1} the proof is complete. \QED

\textbf{Proof of Theorem \ref{th2.3}}
Let $\bigl\{\mathbf{Z}_{n,0}=\bigl(Z_{n,0}^{(1)},Z_{n,0}^{(2)}\bigr),n\geq1 \bigr\}$ be a sequence of 2-dimensions standard Gaussian \ee{random vectors (with mean-zero and unit-variance)} and
\BQNY
corr\left(Z_{n,0}^{(1)},Z_{n,0}^{(2)}\right)=
\frac{\tau_{12}(n)}{\sqrt{\tau_{11}(n)\tau_{22}(n)}}, \quad \tau_{ij}(n)=\frac{\tau_{ij}}{\lne} , \quad i,j \in \{1,2\}.
\EQNY
Further, let $\bigl\{\mathbf{Z}_{n,k}=\bigl(Z_{n,k}^{(1)},Z_{n,k}^{(2)}\bigr), 1\leq k\leq n, n\geq 1 \bigr\}$
denote a triangular array of independent standard Gaussian \ee{random vectors} such that
$corr\left(Z_{n,k}^{(i)},Z_{n,l}^{(j)}\right)=0$ for $1\leq k\neq l\leq n$, $i,j\in\{1,2\}$, and
\begin{eqnarray*}
corr\left(Z_{n,k}^{(1)},Z_{n,k}^{(2)}\right)=
\frac{\rho_0(n)-\tau_{12}(n)}{\sqrt{(1-\tau_{11}(n))(1-\tau_{22}(n))}}=:\tilde{\rho}_0(n).
\end{eqnarray*}
Suppose that $\mathbf{Z}_{n,0}$ is independent of $\{\mathbf{Z}_{n,k}, 1\leq k\leq n\}$, and define
$$\mathbf{Y}_{n,k}=\left(Y_{n,k}^{(1)},Y_{n,k}^{(2)}\right)=
\left(\tau_{11}^{\frac{1}{2}}(n)Z_{n,0}^{(1)}+(1-\tau_{11}(n))^{\frac{1}{2}}Z_{n,k}^{(1)},
\tau_{22}^{\frac{1}{2}}(n)Z_{n,0}^{(2)}+(1-\tau_{22}(n))^{\frac{1}{2}}Z_{n,k}^{(2)}\right)$$
for $1\leq k\leq n$, $n\geq1$. It follows that $\mathbf{Y}_{n,k}$ \ee{are standard Gaussian random vectors}
and $corr\left(Y_{n,k}^{(1)},Y_{n,k}^{(2)}\right)=\rho_{0}(n)$, $corr\left(Y_{n,k}^{(i)},Y_{n,l}^{(j)}\right)=\tau_{ij}(n)$
for $1\leq k\neq l\leq n$, $i,j \in \{1,2\}$.
For notational simplicity, in the sequel denote
$$\left(\widehat{M}_n^{(1)},\widehat{M}_n^{(2)}\right):=\left(\max_{1\leq k\leq n}Y_{n,k}^{(1)},\max_{1\leq k\leq n}Y_{n,k}^{(2)}\right),
\quad
\left(\widetilde{M}_n^{(1)},\widetilde{M}_n^{(2)}\right):=\left(\max_{1\leq k\leq n}Z_{n,k}^{(1)},\max_{1\leq k\leq n}Z_{n,k}^{(2)}\right).$$
\hh{Since $\cw{b_n\sim \sqrt{2\ln n}}$, condition \eqref{eq1.1} implies} 
\COM{
\begin{eqnarray*}
\frac{u_n(x)-\tilde{\rho}_0(n)u_n(t)}{\sqrt{1-\tilde{\rho}_0^2(n)}}=\frac{b_n\sqrt{1-\tilde{\rho}_0(n)}}{\sqrt{1+\tilde{\rho}_0(n)}}+
\frac{(x-t)a_n}{\sqrt{1-\tilde{\rho}_0^2(n)}}+\frac{\sqrt{1-\tilde{\rho}_0(n)}a_nt}{\sqrt{1+\tilde{\rho}_0(n)}}
\end{eqnarray*}
and 
}
\begin{eqnarray*}
\lim_{n \to \infty}\frac{b_n^2(1-\tilde{\rho}_0(n))}{1+\tilde{\rho}_0(n)}&=&
 \widetilde{\lambda}^2,
\end{eqnarray*}
where $ \widetilde{\lambda}=\sqrt{\lambda^2+\tilde{\tau}}$, $\tilde{\tau}=\tau_{12}-\frac{1}{2}(\tau_{11}+\tau_{22})$
and $\lambda^2\cz{\geq}-\tilde{\tau}$. \hh{Consequently, in view of H\"usler and Reiss (1989)}
$$\lim_{n\to \infty}\sup_{x,y\in \R}
\left|\pk{\widetilde{M}_n^{(1)}\leq u_n(x),\widetilde{M}_n^{(2)}\leq u_n(y)}-H_{ \widetilde{\lambda}}(x,y)\right|=0.$$
\cz{
Particularly, for $\tilde{\lambda}=0$ case, i.e., $\lambda^2=-\tilde{\tau}$
we have $1-\rho_0(n)\sim\frac{1}{2}(\tau_{11}(n)+\tau_{22}(n))-\tau_{12}(n).$
Thus $$\lim_{n \to \infty}\tilde{\rho}_0(n)=\lim_{n \to \infty}\frac{1-\frac{1}{2}(\tau_{11}(n)+\tau_{22}(n))}{\sqrt{(1-\tau_{11}(n))(1-\tau_{22}(n))}}= 1,$$
i.e., the asymptotic complete dependence of the components of $\mathbf{Z}_{n,k}$.
For $\tilde{\lambda}=\infty$ case, i.e., $\lambda=\infty$, we have $\limit{n}\rho_0(n)= 0$. Consequently,
$$\lim_{n \to \infty}\tilde{\rho}_0(n)=\lim_{n \to \infty}\frac{-\tau_{12}(n)}{\sqrt{(1-\tau_{11}(n))(1-\tau_{22}(n))}}=0,$$
and thus the asymptotic independence of the components of $\mathbf{Z}_{n,k}$ follows.
}
For all $x,y \in \R$
\begin{eqnarray*}
&&\lim_{n\to \infty}\P\left(\widehat{M}^{(1)}_n\leq u_n(x),\widehat{M}^{(2)}_n\leq u_n(y)\right)\\
&=&\lim_{n\to \infty}\P\left(\tau_{11}^{\frac{1}{2}}(n)Z_{n,0}^{(1)}+(1-\tau_{11}(n))^{\frac{1}{2}}\widetilde{M}^{(1)}_n\leq u_n(x),
\tau_{22}^{\frac{1}{2}}(n)Z_{n,0}^{(2)}+(1-\tau_{22}(n))^{\frac{1}{2}}\widetilde{M}^{(2)}_n\leq u_n(y)\right)\\
&=&\lim_{n\to \infty}\int^{+\infty}_{-\infty}\int^{+\infty}_{-\infty}
\P\left(\widetilde{M}^{(1)}_n\leq (u_n(x)-\tau_{11}^{\frac{1}{2}}(n)z_1)(1-\tau_{11}(n))^{-\frac{1}{2}},\right.\\
&&\qquad \qquad \qquad \qquad \qquad \quad
\left. \widetilde{M}^{(2)}_n\leq (u_n(y)-\tau_{22}^{\frac{1}{2}}(n)z_2)(1-\tau_{22}(n))^{-\frac{1}{2}}\right)\varphi_n(z_1,z_2)dz_1dz_2,
\end{eqnarray*}
 \cw{where $\varphi_n(z_1,z_2)$ is the joint \cH{probability} density of Gaussian vector $\mathbf{Z}_{n,0}$.}
Using \eqref{eq2.2} and $\tau_{ii}(n)=\tau_{ii}/ \lne $, for $i\in \{1,2\}$, let $v_1=x$, $v_2=y$, we have
\BQNY
(u_n(v_i)-\tau_{ii}^{\frac{1}{2}}(n)z_i)(1-\tau_{ii}(n))^{-\frac{1}{2}}=u_n(v_i+\tau_{ii}-\sqrt{2 \tau_{ii}}z_i).
\EQNY
Hence by the dominated convergence theorem, for all $x,y \in \R$ we have
\BQNY
\lim_{n\to \infty}\P\left(\widehat{M}^{(1)}_n\leq u_n(x),\widehat{M}^{(2)}_n\leq u_n(y)\right)
=\EE{H_{ \widetilde{\lambda}}\left(x+\tau_{11}-\sqrt{2\tau_{11}}Z,y+\tau_{22}-\sqrt{2\tau_{22}}W\right)},
\EQNY
where $\left(Z,W\right)$ is a standard Gaussian vector with correlation $\tau_{12}/\sqrt{\tau_{11}\tau_{22}}$.\\
\hh{Applying} \aA{Berman's} Normal Comparison Lemma and  \eqref{L2}, we get 
\begin{eqnarray*}
&&\left|\P\left(M^{(1)}_n\leq u_n(x),M^{(2)}_n\leq u_n(y)\right)
-\P\left(\widehat{M}^{(1)}_n\leq u_n(x),\widehat{M}^{(2)}_n\leq u_n(y)\right)\right|\\
\COM{
&\leq& \CC{\frac{1}{4} \sum_{1\leq k<l \leq n}  |\Cov(X_{n,k}^{(1)},X_{n,l}^{(1)})-\Cov(Y_{n,k}^{(1)},Y_{n,l}^{(1)})|\exp\left(-\frac{u_n^2(x)}{1+\varrho_{11}(k,l,n)}\right)}\\
&&\CC{+\frac{1}{4}\sum_{1\leq k\neq l\leq n}|\Cov(X_{n,k}^{(1)},X_{n,l}^{(2)})-\Cov(Y_{n,k}^{(1)},Y_{n,l}^{(2)})|\exp\left(-\frac{u_n^2(x)+u_n^2(y)}{2(1+\varrho_{12}(k,l,n))}\right)}\\
&&\CC{+\frac{1}{4}\sum_{1\leq k<l \leq n}|\Cov(X_{n,k}^{(2)},X_{n,l}^{(2)})-\Cov(Y_{n,k}^{(2)},Y_{n,l}^{(2)})|\exp\left(-\frac{u_n^2(y)}{1+\varrho_{22}(k,l,n)}\right)}\\
}
&\leq& \frac{1}{4} n\sum^{n-1}_{k=1}  |\rho_{11}(k,n)-\tau_{11}(n)|\exp\left(-\frac{u_n^2(x)}{1+\varrho_{11}(k,n)}\right)\\
&&+\frac{1}{2}n\sum^{n-1}_{k=1}|\rho_{12}(k,n)-\tau_{12}(n)|\exp\left(-\frac{u_n^2(x)+u_n^2(y)}{2(1+\varrho_{12}(k,n))}\right)\\
&&+\frac{1}{4}n\sum^{n-1}_{k=1}|\rho_{22}(k,n)-\tau_{22}(n)|\exp\left(-\frac{u_n^2(y)}{1+\varrho_{22}(k,n)}\right)\\
&\to&0,\qquad \hh{n\to \IF},
\end{eqnarray*}
thus the claim follows. \QED

\textbf{Proof of Theorem \ref{th1.2}}.
Let $\left\{\left(\hat{X}_{n,k}^{(1)},\hat{X}_{n,k}^{(2)}\right), 1\leq k\leq n, n\geq1\right\}$ be
iid Gaussian triangular array defined in the proof of Theorem \ref{th2.1}, and
set $\left(\mathfrak{m}_n^{(1)},\mathfrak{m}_n^{(2)}\right):=\left(\min_{1\leq k \leq n}\hat{X}_{n,k}^{(1)}, \min_{1\leq k \leq n}\hat{X}_{n,k}^{(2)}\right)$. \aA{We write next}
\BQNY
\lefteqn{n\left(1-\pk{-u_n(y_1)<\hat{X}^{(1)}_{1,1}\leq u_n(x_1),-u_n(y_2)< \hat{X}^{(2)}_{1,1}\leq u_n(x_2)}\right)}\\
&
=&nP_1(n,x_1,x_2)+nP_2(n,y_1,y_2)-nP_3(n,x_1,y_2)-nP_4(n,y_1,x_2),
\EQNY
\COM{
\BQNY
&&\pk{-u_n(y_1)< \mathfrak{m}_n^{(1)}\leq \mathfrak{M}_n^{(1)}\leq u_n(x_1),
-u_n(y_2)< \mathfrak{m}_n^{(2)}\leq \mathfrak{M}_n^{(2)}\leq u_n(x_2)}\\
&=&\EXP{-n\left(1-\pk{-u_n(y_1)<\hat{X}_1\leq u_n(x_1),-u_n(y_2)< \hat{X}_2\leq u_n(x_2)}\right)}+o(1)\\
&=&\EXP{-nP_1(n,x_1,x_2)-nP_2(n,y_1,y_2)+nP_3(n,x_1,y_2)+nP_4(n,y_1,x_2)}+o(1),
\EQNY}
where
\BQNY
&&P_1(n,x_1,x_2):=\pk{\hat{X}^{(1)}_{1,1}>u_n(x_1)}+\pk{\hat{X}^{(2)}_{1,1}>u_n(x_2)}-
\pk{\hat{X}^{(1)}_{1,1}>u_n(x_1),\hat{X}^{(2)}_{1,1}>u_n(x_2)},\\
&&P_2(n,y_1,y_2):=\pk{\hat{X}^{(1)}_{1,1}\leq-u_n(y_1)}+\pk{\hat{X}^{(2)}_{1,1}\leq-u_n(y_2)}-
\pk{\hat{X}^{(1)}_{1,1}\leq-u_n(y_1),\hat{X}^{(2)}_{1,1}\leq-u_n(y_2)},\\
&&P_3(n,x_1,y_2):=\pk{\hat{X}^{(1)}_{1,1}>u_n(x_1),\hat{X}^{(2)}_{1,1}\leq-u_n(y_2)},\\
&&P_4(n,y_1,x_2):=\pk{\hat{X}^{(1)}_{1,1}\leq-u_n(y_1),\hat{X}^{(2)}_{1,1}>u_n(x_2)}.
\EQNY
Using Theorem 1 in H\"{u}sler and Reiss (1989), we
have
\BQNY
&&\lim_{n \to \infty}nP_1(n,x_1,x_2)
=\Phi\left(\lambda+\frac{x_1-x_2}{2\lambda}\right)\exp(-x_2)+
\Phi\left(\lambda+\frac{x_2-x_1}{2\lambda}\right)\exp(-x_1)=:D_1,
\EQNY
and since $(-\hat{X}^{(1)}_{1,1},-\hat{X}^{(2)}_{1,1})\stackrel{d}{=}(\hat{X}^{(1)}_{1,1},\hat{X}^{(2)}_{1,1})$
\BQNY
&&\lim_{n \to \infty}nP_2(n,y_1,y_2)
=\Phi\left(\lambda+\frac{y_1-y_2}{2\lambda}\right)\exp(-y_2)+
\Phi\left(\lambda+\frac{y_2-y_1}{2\lambda}\right)\exp(-y_1)=:D_2.
\EQNY
In view of \eqref{eq2.2} we have $\lim_{n \to \infty}\frac{n a_n}{\sqrt{2\pi}} \EXP{-\frac{b_n^2}{2}}=1$, then \aA{with
$\CCC{\varphi}= \Phi'$}
\BQNY
nP_3(n,x_1,y_2)
&=&n a_n \int_{x_1}^{\infty}\Phi\left(\frac{-u_n(y_2)-\rho_0(n)u_n(t)}{\sqrt{1-\rho_0^2(n)}}\right)\varphi(u_n(t))dt\\
&\sim&\int_{x_1}^{\infty}\Phi\left(\frac{-u_n(y_2)-\rho_0(n)u_n(t)}{\sqrt{1-\rho_0^2(n)}}\right)\EXP{-a_nb_nt-\frac{a_n^2t^2}{2}}dt, \quad \aA{n\to \IF}.
\EQNY
By \eqref{eq1.1} and \eqref{eq2.2}
\BQNY
\frac{-u_n(y_2)-\rho_0(n)u_n(t)}{\sqrt{1-\rho_0^2(n)}}=
-\frac{b_n\sqrt{1+\rho_0(n)}}{\sqrt{1-\rho_0(n)}}-\frac{a_n(y_2+t)}{\sqrt{1-\rho_0^2(n)}}+\frac{a_nt\sqrt{1-\rho_0(n)}}{\sqrt{1+\rho_0(n)}}\to -\infty
\EQNY
as $n\to \infty$, hence  $\lim_{n \to \infty}nP_3(n,x_1,y_2)=0,$  and by similar arguments $\lim_{n \to \infty}nP_4(n,y_1,x_2)=0.$
Consequently, for all $x_1,x_2,y_1,y_2\inr$
\BQNY
&&\lim_{n \to \infty}\pk{-u_n(y_1)< \mathfrak{m}_n^{(1)}\leq \mathfrak{M}_n^{(1)}\leq u_n(x_1),-u_n(y_2)< \mathfrak{m}_n^{(2)}\leq \mathfrak{M}_n^{(2)}\leq u_n(x_2)}\\
&=&\lim_{n \to \infty}\left(\CCC{\pk{-u_n(y_1)<\hat{X}^{(1)}_{1,1}\leq u_n(x_1),-u_n(y_2)< \hat{X}^{(2)}_{1,1}\leq u_n(x_2)}}\right)^n\\
&=&\lim_{n \to \infty}\left(1-\frac{D_1+D_2}{n}+o\left(\frac{1}{n}\right)\right)^n\\
&=&\EXP{-D_1-D_2}=H_{\lambda}(x_1,x_2)H_{\lambda}(y_1,y_2).
\EQNY
\COM{
\BQNY
\lim_{n \to \infty}\sup_{x_1,x_2,y_1,y_2\inr}
\Biggl \lvert \pk{-u_n(y_1)< \mathfrak{m}_n^{(1)}\leq \mathfrak{M}_n^{(1)}\leq u_n(x_1),-u_n(y_2)< \mathfrak{m}_n^{(2)}\leq \mathfrak{M}_n^{(2)}\leq u_n(x_2)}
\notag\\
&&\hspace{-11 cm} - H_{\lambda}(x_1,x_2)H_{\lambda}(y_1,y_2)\Biggr \lvert =0.
\EQNY}
By Berman's Normal Comparison Lemma 
and \eqref{L1}, for all $x_1,x_2,y_1,y_2 \in \R$
with
$$\omega_n:=\min(|u_n(x_1)|,|u_n(x_2)|,|u_n(y_1)|,|u_n(y_2)|)$$
 we have
\begin{eqnarray*}
&&\left|\pk{-u_n(y_1)< m_n^{(1)}\leq M_n^{(1)}\leq u_n(x_1),-u_n(y_2)< m_n^{(2)}\leq M_n^{(2)}\leq u_n(x_2)}\right.\\
&&\quad -\left.\pk{-u_n(y_1)<  \mathfrak{m}_n^{(1)}\leq \mathfrak{M}_n^{(1)}\leq u_n(x_1),-u_n(y_2)< \mathfrak{m}_n^{(2)}\leq \mathfrak{M}_n^{(2)}\leq u_n(x_2)}\right|\\
\COM{&\leq&\CC{\sum_{1\leq k < l \leq n}|\Cov(X_{n,k}^{(1)},X_{n,l}^{(1)})|\exp\left(-\frac{\omega^2_n}{1+|\Cov(X_{n,k}^{(1)},X_{n,l}^{(1)})|}\right)}\\
&&\CC{+
\sum_{1\leq k \neq l \leq n}|\Cov(X_{n,k}^{(1)},X_{n,l}^{(2)})|\exp\left(-\frac{\omega^2_n}{1+|\Cov(X_{n,k}^{(1)},X_{n,l}^{(2)})|}\right)}\\
&&+\CC{\sum_{1\leq k < l \leq n}|\Cov(X_{n,k}^{(2)},X_{n,l}^{(2)})|\exp\left(-\frac{\omega^2_n}{1+|\Cov(X_{n,k}^{(2)},X_{n,l}^{(2)})|}\right)}\\
}
&\leq&n\sum_{k=1}^{n-1}|\rho_{11}(k,n)|\exp\left(-\frac{\omega^2_n}{1+|\rho_{11}(k,n)|}\right)+
2n\sum_{k=1}^{n-1}|\rho_{12}(k,n)|\exp\left(-\frac{\omega^2_n}{1+|\rho_{12}(k,n)|}\right)\\
&&+n\sum_{k=1}^{n-1}|\rho_{22}(k,n)|\exp\left(-\frac{\omega^2_n}{1+|\rho_{22}(k,n)|}\right)\\
&\to & 0,  \quad n \to \infty,
\end{eqnarray*}
and hence the claim follows. \QED

\textbf{Proof of Theorem \ref{th2.2}.}
 \hh{In the light} of Theorem \ref{th2.1}, \hh{\eqref{eq2.4} follows if we show that} for any $x,y\inr$
\begin{eqnarray}\label{eq3.4}
\lim_{n \to \infty}\frac{1}{\ln n}\sum^n_{k=1}\frac{1}{k}\left(\I\left(M^{(1)}_k\leq u_k(x), M^{(2)}_k\leq u_k(y)\right)
-\P\left(M^{(1)}_k\leq u_k(x), M^{(2)}_k\leq u_k(y)\right)\right) = 0
\end{eqnarray}
holds almost surely. 
\hh{By} Lemma 3.1 in Cs\'{a}ki and Gonchigdanzan (2002), in order to prove \eqref{eq3.4},
it suffices to show that \cH{for some} $\epsilon>0$ and \cw{some positive constant $c$}
$$\Var\left(\sum^n_{k=1}\frac{1}{k}\I\left(M^{(1)}_k\leq u_k(x), M^{(2)}_k\leq u_k(y)\right)\right)\leq
c(\ln n)^2(\ln\ln n)^{-(1+\epsilon)}.
$$
Straightforward calculations yield
\begin{eqnarray*}
&&\Var\left(\sum^n_{k=1}\frac{1}{k}\I\left(M^{(1)}_k\leq u_k(x), M^{(2)}_k\leq u_k(y)\right)\right)\\
&=&\sum^n_{k=1}\frac{1}{k^2}\Var\left(\I\left(M^{(1)}_k\leq u_k(x), M^{(2)}_k\leq u_k(y)\right)\right)\\
&&+2\sum^n_{1\leq k< l\leq n}\frac{1}{kl}\Cov\left(\I\left(M^{(1)}_k\leq u_k(x), M^{(2)}_k\leq u_k(y)\right),
\I\left(M^{(1)}_l\leq u_l(x), M^{(2)}_l\leq u_l(y)\right)\right)\\
&=:&\Sigma_1+\Sigma_2,
\end{eqnarray*}
where $\Sigma_1<\infty$. Thus, it remains only to estimate $\Sigma_2$. Write next
\begin{eqnarray*}
&&\Cov\left(\I\left(M^{(1)}_k\leq u_k(x), M^{(2)}_k\leq u_k(y)\right),
\I\left(M^{(1)}_l\leq u_l(x), M^{(2)}_l\leq u_l(y)\right)\right)\\
&\leq&2\EE{\left|\I\left(M^{(1)}_l\leq u_l(x),M^{(2)}_l\leq u_l(y)\right)-\I\left(M^{(1)}_{l,k}\leq u_l(x),M^{(2)}_{l,k}\leq u_l(y)\right)\right|}\\
&&+\left|\Cov\left(\I\left(M^{(1)}_k\leq u_k(x),M^{(2)}_k\leq u_k(y)\right),\I\left(M^{(1)}_{l,k}\leq u_l(x),M^{(2)}_{l,k}\leq u_l(y)\right)\right)\right|\\
&=:&2E_1+E_2.
\end{eqnarray*}
For $1\leq k<l\leq n$, we have
\begin{eqnarray*}
E_1
&\leq&\left|\P\left(M^{(1)}_{l,k}\leq u_l(x),M^{(2)}_{l,k}\leq u_l(y)\right)
-\P\left(\mathfrak{M}^{(1)}_{l,k}\leq u_l(x),\mathfrak{M}^{(2)}_{l,k}\leq u_l(y)\right)\right|\\
&&+\left|\P\left(M^{(1)}_{l}\leq u_l(x),M^{(2)}_{l}\leq u_l(y)\right)
-\P\left(\mathfrak{M}^{(1)}_{l}\leq u_l(x),\mathfrak{M}^{(2)}_{l}\leq u_l(y)\right)\right|\\
&&+\left|\P\left(\mathfrak{M}^{(1)}_{l,k}\leq u_l(x),\mathfrak{M}^{(2)}_{l,k}\leq u_l(y)\right)
-\P\left(\mathfrak{M}^{(1)}_{l}\leq u_l(x),\mathfrak{M}^{(2)}_{l}\leq u_l(y)\right)\right|\\
&=:&P_1+P_2+P_3,
\end{eqnarray*}
\hh{where}
$$\left(\mathfrak{M}_{l,k}^{(1)}, \mathfrak{M}_{l,k}^{(2)}\right):
=\left(\max_{k+1\leq t\leq l }\hat{X}^{(1)}_{l,t}, \max_{k+1\leq
t\leq l }\hat{X}^{(2)}_{l,t}\right), \quad
\left(\mathfrak{M}_{l}^{(1)},
\mathfrak{M}_{l}^{(2)}\right):=\left(\mathfrak{M}_{l,0}^{(1)},
\mathfrak{M}_{l,0}^{(2)}\right).$$ \hh{By} \CCC{Berman's} Normal Comparison
Lemma and \eqref{eq3.1}, \bC{for all $l\le n$} we obtain
\BQN\label{add eq2}
\aA{P_1+P_2}\leq c(\ln \ln l)^{-(1+\epsilon)},
\EQN \cH{with $c$ some positive constant.}
\hh{The following} inequality $z^{l-k}-z^{l}\leq\frac{k}{l}$ \hh{valid for} $0\leq z \leq1$ \hh{yields further }
\BQN\label{add eq3}
P_3=H^{l-k}(u_l(x),u_l(y))-H^{l}(u_l(x),u_l(y))\leq \frac{k}{l},
\EQN
where $H$ \cH{is} the \cH{distribution function} of $\left(\hat{X}^{(1)}_{11},\hat{X}^{(2)}_{11}\right)$.
\aA{Again, as above by} \eqref{eq3.2}, we \hh{have with 
\cZ{$c$ some} positive constant}
\begin{eqnarray}\label{add eq4}
E_2
&\leq&\left|\P\left(M^{(1)}_k\leq u_k(x),M^{(2)}_k\leq u_k(y),M^{(1)}_{l,k}\leq u_l(x),M^{(2)}_{l,k}\leq u_l(y)\right)\right.\nonumber\\
&&-\left.\P\left(M^{(1)}_k\leq u_k(x),M^{(2)}_k\leq u_k(y)\right)
\P\left(M^{(1)}_{l,k}\leq u_l(x),M^{(2)}_{l,k}\leq u_l(y)\right)\right|\nonumber\\
\COM{&\leq&
\CC{\frac{1}{4} \sum_{i,j=1,2} \sum_{1\leq s\leq k \atop k+1\leq t\leq l} |\Cov(X_{k,s}^{(i)},X_{l,t}^{(j)})|
\exp\left(-\frac{u_k^2(z_i)+u_l^2(z_j)}{2(1+|\Cov(X_{k,s}^{(i)},X_{l,t}^{(j)})|)}\right)}\nonumber\\
}
&\leq&
\hh{\frac{1}{4} \sum_{i,j=1,2} \sum_{1\leq s\leq k \atop k+1\leq t\leq l} |\gamma_{ij}(s,t,k,l)|
\exp\left(-\frac{u_k^2(z_i)+u_l^2(z_j)}{2(1+|\gamma_{ij}(s,t,k,l)|)}\right)}\nonumber\\
\COM{&\leq&\mathbb{C}\left[
\sum_{1\leq s\leq k \atop k+1\leq t\leq l}|\gamma_{11}(s,t,k,l)|
\exp\left(-\frac{u_k^2(x)+u_l^2(x)}{2(1+|\gamma_{11}(s,t,k,l)|)}\right)
+
\sum_{1\leq s\leq k \atop k+1\leq t\leq l}|\gamma_{12}(s,t,k,l)|
\exp\left(-\frac{u_k^2(x)+u_l^2(y)}{2(1+|\gamma_{11}(s,t,k,l)|)}\right)\right.\nonumber\\
&&\left.+
\sum_{1\leq s\leq k \atop k+1\leq t\leq l}|\gamma_{21}(s,t,k,l)|
\exp\left(-\frac{u_k^2(y)+u_l^2(x)}{2(1+|\gamma_{11}(s,t,k,l)|)}\right)
+
\sum_{1\leq s\leq k \atop k+1\leq t\leq l}|\gamma_{22}(s,t,k,l)|
\exp\left(-\frac{u_k^2(y)+u_l^2(y)}{2(1+|\gamma_{11}(s,t,k,l)|)}\right)\right]
\nonumber\\
}
&\leq& \frac{1}{4}\hh{k\sum_{i,j=1,2}\sum^l_{s=1} |\gamma_{ij}(s,k,l)|\exp\left(-\frac{u_k^2(z_i)+u_l^2(z_j)}{2(1+|\gamma_{ij}(s,k,l)|)}\right)}\nonumber\\
\COM{&\leq& \mathbb{C}k\sum^l_{s=1}|\gamma_{11}(s,k,l)|exp\left(-\frac{u_k^2(x)+u_l^2(x)}{2(1+|\gamma_{11}(s,k,l)|)}\right)
+\mathbb{C}k\sum^l_{s=1}|\gamma_{12}(s,k,l)|exp\left(-\frac{u_k^2(x)+u_l^2(y)}{2(1+|\gamma_{12}(s,k,l)|)}\right)\nonumber\\
&&+\mathbb{C}k\sum^l_{s=1}|\gamma_{21}(s,k,l)|exp\left(-\frac{u_k^2(y)+u_l^2(x)}{2(1+|\gamma_{21}(s,k,l)|)}\right)
+\mathbb{C}k\sum^l_{s=1}|\gamma_{22}(s,k,l)|exp\left(-\frac{u_k^2(y)+u_l^2(y)}{2(1+|\gamma_{22}(s,k,l)|)}\right)\nonumber\\
}
&\leq&c(\ln \ln l)^{-(1+\epsilon)},
\end{eqnarray}
where $z_1=x$ and $z_2=y$.
Combining \eqref{add eq2}, \eqref{add eq3} and \eqref{add eq4},
for $1\leq k< l\leq n$ we have
$$\Cov\left(\I\left(M^{(1)}_k\leq u_k(x), M^{(2)}_k\leq u_k(y)\right),
\I\left(M^{(1)}_l\leq u_l(x), M^{(2)}_l\leq u_l(y)\right)\right)\leq
c\left(\frac{k}{l}+(\ln \ln l)^{-(1+\epsilon)}\right),$$ hence
$\Sigma_2\leq c(\ln n)^2(\ln\ln n)^{-(1+\epsilon)}$ \cZ{with some
positive constant $c$}. \bC{The proof of \eqref{eq2.4} is
complete.}

\bC{By using Theorem \ref{th1.2} and arguments similar to the proof
of \eqref{eq2.4}, we can show that \eqref{eq2.7} holds. The details
are omitted here.}
 \QED

\cH{\textbf{Acknowledgments}}: Z. Weng kindly acknowledges support by the Swiss National Science Foundation Grant 200021-134785. E. Hashorva kindly acknowledges support by the Swiss National Science Foundation Grant 200021-140633/1.

\end{document}